\documentclass[11pt,a4paper]{article}

\usepackage[left=2.4cm, top=2.45cm,bottom=2.45cm,right=2.4cm]{geometry}
\usepackage{amsmath}
\usepackage{amssymb,amsthm,graphicx,xcolor,tikz,pgfplots,bm}
\usepackage{bbm}

\usepackage[british]{babel}
\usepackage{enumerate}

\usepackage[
bookmarks=true,
bookmarksnumbered=true,
bookmarksopen=true,
unicode=true,
pdftoolbar=true,
pdfmenubar=true,
pdffitwindow=false,
pdfstartview={FitH},
pdftitle={},
pdfauthor={},
pdfsubject={},
pdfcreator={},
pdfproducer={},
pdfkeywords={},
pdfnewwindow=true,
colorlinks=true,
linkcolor=black,
citecolor=black,
filecolor=black,
urlcolor=black]{hyperref}

\makeatletter
\def\@cite#1#2{[\textbf{#1}\if@tempswa , #2\fi]}	
\def\@biblabel#1{[#1]}								
\makeatother

\newtheorem{theorem}{Theorem}
\newtheorem{proposition}[theorem]{Proposition}
\newtheorem{lemma}[theorem]{Lemma}

\theoremstyle{definition}

\newtheorem {remark}[theorem]{Remark}

\newcommand{\EE}{\mathbb{E}}
\newcommand{\HH}{\mathbb{H}}

\newcommand{\PP}{\mathbb{P}}
\newcommand{\RR}{\mathbb{R}}
\renewcommand{\SS}{\mathbb{S}}
\renewcommand{\AA}{\mathbb{A}}

\newcommand{\cH}{\mathcal{H}}

\DeclareMathOperator{\Var}{Var}
\DeclareMathOperator{\Wass}{Wass}
\DeclareMathOperator{\Kol}{Kol}
\DeclareMathOperator{\cum}{cum}

\DeclareMathOperator{\asinh}{arcsinh}
\newcommand{\bo}{\mathbf{o}}

\newcommand{\eps}{\varepsilon}
\newcommand{\dint}{\mathrm{d}}

\newcommand\tinyonehalf{\scalebox{0.65}{$\frac{1}{2}$}}

\begin{document}
	
	\title{\bfseries A quantitative central limit theorem for\\ Poisson horospheres in high dimensions} \date{}
	
	\author{%
		Zakhar Kabluchko\footnotemark[1]%
		\and Daniel Rosen\footnotemark[2]%
		\and Christoph Th\"ale\footnotemark[3]%
	}
	
	\date{}
	\renewcommand{\thefootnote}{\fnsymbol{footnote}}
	\footnotetext[1]{%
		University of M\"unster, Germany. Email: zakhar.kabluchko@uni-muenster.de
	}
	
	\footnotetext[2]{%
		Ruhr University Bochum, Germany. Email: daniel.rosen@rub.de
	}	
	
	\footnotetext[3]{%
		Ruhr University Bochum, Germany. Email: christoph.thaele@rub.de
	}

	\maketitle
	
	\begin{abstract}
		\noindent Consider a stationary Poisson process of horospheres in a $d$-dimensional hyperbolic space. In the focus of this note is the total surface area these random horospheres induce in a sequence of balls of growing radius $R$. The main result is a quantitative, non-standard central limit theorem for these random variables as the radius $R$ of the balls and the spatial dimension $d$ tend to infinity simultaneously. \\
		
	\end{abstract}
	
	\section{Introduction and main result}
	
	The study of random geometric systems in non-Euclidean geometries is a recent and fast growing branch of stochastic geometry. We refer to \cite{BenjaminiJonassonSchrammTykesson,BenjaminiPaquettePfeffer,BenjaminiSchramm,BesauRosenThaele,BesauThaele,FountoulakisMuller,FountoulakisYukich,GodlandKabluchkoThaele,HansenMuller,HeroldHugThaele,OwadaYogesh} for selected works on hyperbolic random geometric graphs, random tessellations and random polytopes.
	
	In this note we address an interesting generalization of the Poisson hyperplane process to hyperbolic geometry. The study of Euclidean Poisson hyperplanes is by now classical \cite{Heinrich,LPST,RS,SchneiderWeil} and was extended in \cite{HeroldHugThaele} to hyperbolic space, where {Poisson processes} of totally geodesic hypersurfaces {are} studied, see also \cite{SantaloYanez} for mean values in the planar case. Even more recently, in \cite{KRT} it was observed that this model fits into a one-parameter family of so-called \emph{Poisson $\lambda$-geodesic hyperplanes}, and the fluctuations of the total hyperbolic surface area of such a process within a sequence of growing balls were examined in detail. The special case we consider here is the \emph{Poisson horosphere process}, {which corresponds to the choice $\lambda = 1$ in \cite{KRT}}. {In that paper it was shown that, in contrast with other choices of $\lambda \in [0,1)$, these fluctuations are Gaussian in all spatial dimensions. The purpose of the present note is to quantify this result and to extend it to the setting of growing dimension.}
	
\begin{figure}[t]
		\centering
		\includegraphics[width=.4\textwidth]{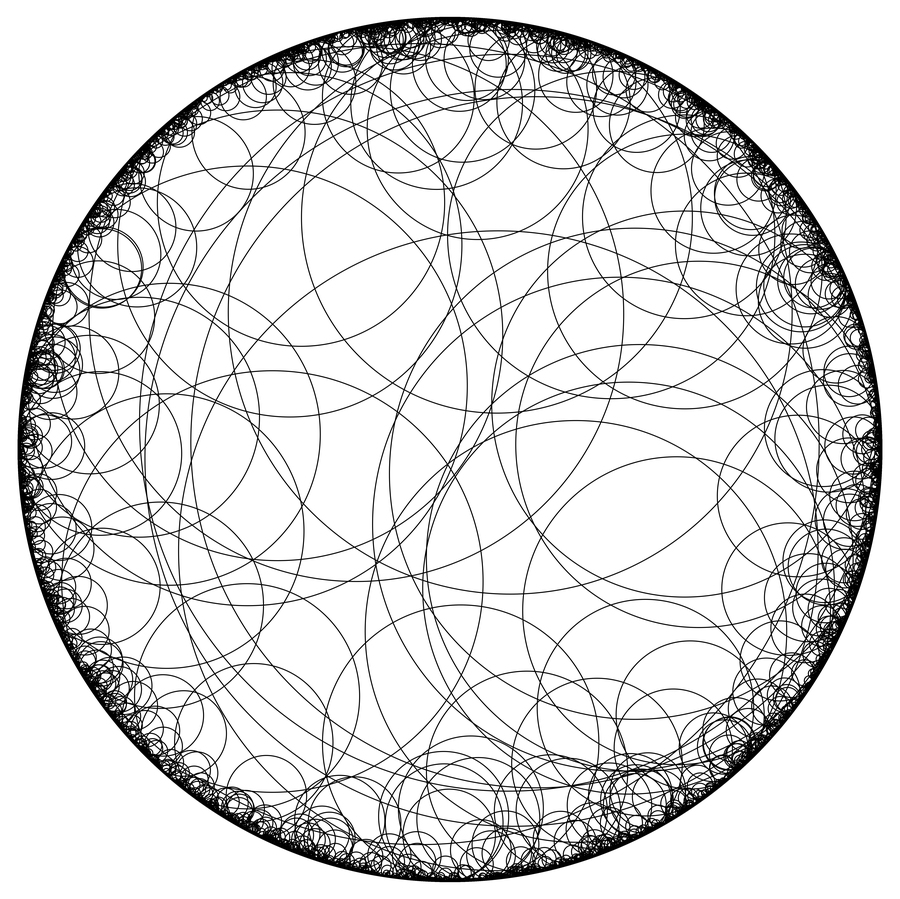}
		\caption{Simulation of a Poisson process of horospheres in the Poincar\'e disc model for the hyperbolic plane.}\label{fig:Poisson_horospheres}
	\end{figure}

	Let us recall some definitions; for more details we refer the reader to \cite{KRT} and the references cited therein. A \emph{horosphere} in a $d$-dimensional hyperbolic space $\HH^d$ is, intuitively speaking, a sphere of infinite radius. More formally, it is a complete totally umbilic hypersurface of constant normal curvature $1$. For concreteness, in the Poincar\'e ball model of hyperbolic space, horospheres are realized as Euclidean spheres tangent to the boundary, see Figure \ref{fig:Poisson_horospheres}.  We denote by $\cH$ the space of all horospheres in $\HH^d$. This space admits a transitive action by the group of hyperbolic isometries and an invariant measure for this action, which is unique up to a multiplicative constant and will be denoted by $\Lambda$, see \cite{GNS,Sol}.
	
	Now, let $\eta_d$ be a Poisson process on $\cH$ with intensity measure $\Lambda$, see Figure \ref{fig:Poisson_horospheres} for a simulation in the case $d=2$. 
	For $R > 0$, we consider the total surface area
	\begin{equation*}
		S_{R,d} := \sum_{H \in \eta_d} \cH^{d-1}(H \cap B^d_R)
	\end{equation*}
	of $\eta_d$ within a hyperbolic ball $B^d_R$ around an arbitrary but fixed point in $\HH^d$ {and having} hyperbolic radius $R>0$. Here, $\cH^{d-1}$ stands for the $(d-1)$-dimensional Hausdorff measure with respect to the hyperbolic metric. In \cite{KRT} it was proven that, for a fixed spatial dimension $d$, the centred and normalized surface area satisfies a non-standard central limit theorem. Namely, it converges in distribution, as $R\to\infty$ and after centering and normalizing by the standard deviation, to a Gaussian random variable of variance ${1 \over 2}$. The main result of the present note extends this in two directions: first, we provide estimates on the rate of convergence. Second, our bounds depend explicitly on the dimension, providing central limit theorems for Poisson horospheres in increasing spatial dimensions. To measure the distance between two random variables $X$ and $Y$ we use the {Kolmogorov and Wasserstein metrics, which are defined, respectively, by
	$$
	d_{\Kol}(X,Y) := \sup_{t\in\RR} |\PP(X\leq t)-\PP(Y\leq t)|
	$$
	and
	$$
	d_{\Wass}(X,Y) := \sup_h\big|\EE[h(X)]-\EE[h(Y)]\big|,
	$$
	where the latter supremum is taken over all Lipschitz functions $h:\RR\to\RR$ with Lipschitz constant at most one.}

	{Our main result is the following non-standard quantitative central limit theorem.
	\begin{theorem}\label{thm:CLT-master}
		Let $N_{\tinyonehalf}$ be a centred Gaussian random variable of variance $\frac{1}{2}$. Consider the surface functional $S_{R,d}$, for $d\geq 2$ and $R\geq 1$. Then there exists a universal constant $C > 0$ such that for any choice $\bullet \in \{\Kol,\Wass\}$ one has
		\begin{equation*}
		d_{\bullet}\left(\frac{ S_{R,d} - \EE S_{R,d}}{\sqrt{\Var S_{R,d}}}, N_{\tinyonehalf} \right) \leq  C  \cdot \begin{cases}
		 e^{-R/2} &: R-\log d \leq 1, \\
		   {1 \over \sqrt{d} \,(R - \log d)} + {1 \over d \,\sqrt{R-\log d}} &: R-\log d > 1.
		\end{cases}
		\end{equation*}
	\end{theorem}}

	\begin{remark}\label{rem:main_thm}
	\begin{enumerate}[(i)]
		\item {We first point out that in a fixed spatial dimension $d\geq 2$, Theorem \ref{thm:CLT-master} implies the following bound:
		\begin{align*}
		d_{\bullet}\left(\frac{ S_{R,d} - \EE S_{R,d}}{\sqrt{\Var S_{R,d}}}, N_{\tinyonehalf} \right) & \leq  C_d \, R^{-1/2},
		\end{align*}
		with a constant $C_d$ depending only on $d$.}\\
		We note that the rate of convergence $R^{-1/2}$ is the same in all dimensions. The same convergence rate, in all dimensions, is observed in the central limit theorem for the total surface area of Poisson hyperplanes in Euclidean space (see \cite{LPST}, {but the result is also a special case of \eqref{eq:Wass_bd_Euclidean} below}). Let us remark that this limiting behaviour is in sharp contrast with the cases of $\lambda$-geodesic hyperplanes with $\lambda < 1$ (we recall that horospheres correspond to $\lambda=1$). As described in \cite{KRT}, in those cases the fluctuations of the surface functional are non-Gaussian in every fixed dimension $\geq 4$. The geometric distinction between the two cases is that horospheres are \emph{intrinsically Euclidean}, while the intrinsic geometry of other  $\lambda$-geodesic hyperplanes is hyperbolic. 
		
		\item{ The threshold between the two cases in Theorem \ref{thm:CLT-master} is somewhat arbitrary. The constant $1$ can be replaced by any positive number, at the cost of changing the constant $C$ in the theorem.}
		
		\item We note that the convergence rate in the first case of Theorem \ref{thm:CLT-master} is always worse than $d^{-1/2}$ (and in particular, worse than in the second case). Indeed, by assumption, $e^{R} \leq e \, d$ and hence $e^{-R/2}$ converges to zero slower than $d^{-1/2}$.
		
		\item{ In the high-dimensional regime, that is if $d\to\infty$ and $R=R_d$ is a sequence of radii, we note that Theorem \ref{thm:CLT-master} implies a non-standard central limit theorem for the surface functional as soon as $R_d \xrightarrow[d\to\infty]{}\infty$.} For example, taking the radius as $R=R_d = \alpha \log d$ for $\alpha > 0$, the theorem gives
		\[
		{d_{\bullet}}\left( \frac{S_{R,d} - \EE S_{R,d}  } {\sqrt{\Var S_{R,d} } }, N_{\tinyonehalf} \right)   \leq C \cdot \begin{cases} d^{-\alpha/2} &: \alpha \leq 1 \\ d^{-1/2} (\log d)^{-1} &: \alpha > 1.	
		\end{cases}
		\]		
		It is also natural to ask for sharp conditions on $R_d$ which ensure that the centred and normalized total surface area is asymptotically Gaussian. As noted above, $R_d \to \infty$ is sufficient, but for fixed $R$ our bounds do not yield a central limit theorem for the surface functional. We have to leave this as an open problem.
		
		\item The reader might be interested in a comparison with the Euclidean case, where one considers the total surface area {$S_{R,d,e}$} of a stationary and isotropic Poisson process on the space of hyperplanes in $\RR^d$ within a centred ball of radius $R$. In this situation it holds that
		\begin{equation}\label{eq:Wass_bd_Euclidean}
		{d_{\bullet}}\Big({S_{R,d,e}-\EE S_{R,d,e}\over\sqrt{\Var S_{R,d,e}}},N\Big) \leq C\,d^{1/4}R^{-1/2}
		\end{equation}
		for some absolute constant $C>0$ and where $N$ denotes a standard Gaussian random variable. Since we could not locate this result in the literature, we provide an argument in Section \ref{sec:Euclidean}. {In particular, as $d\to\infty$ we need that $R$ grows faster than $\sqrt{d}$ in order to deduce from \eqref{eq:Wass_bd_Euclidean} a central limit theorem.}
	\end{enumerate}
\end{remark}
	
	\section{Proof of the main results}
	
	Before proving Theorem \ref{thm:CLT-master}, we need to recall some preliminaries. First we need an explicit description of the invariant measure $\Lambda$ on the space $\cH$ of horospheres. We fix an origin $\bo\in \HH^d$ and parametrize an element $H\in \cH$ by the pair $(s, u) \in \RR \times\SS^{d-1}$, where $s\in\RR$ is the signed distance from $H$ to $\bo$ ({with $s > 0$ if $\bo$ lies on the convex side of $H$, and negative otherwise}), and $u\in\SS^{d-1}$ is the unit vector (in the tangent space $T_\bo\HH^d$) along the geodesic passing through $\bo$ and intersecting $H$ orthogonally, {while pointing outside of the convex side}. The invariant measure is then defined by the relation
	\begin{equation}\label{eq:dHlambda}
		\int_{\cH}f(H)\,\Lambda(\dint H) = \int_{\RR}\int_{\SS^{d-1}} f(H(s,u))\,e^{-(d-1)s} \,\dint u\,\dint s,
	\end{equation}
	where $f:\cH\to\RR$ is a non-negative measurable function and $H(s,u)$ stands for the unique element of $\cH$ parametrized by $(s,u)$ as just described. {Here $\dint s$ and $\dint u$ stand for the Lebesgue measure on $\RR$ and the normalized spherical Lebesgue measure on $\SS^{d-1}$, respectively.}
	
	We will also need the following geometric computation of the volume of the intersection $H(s) \cap B_R^d $, where $H(s) \subset \HH^d$ is a horosphere of signed distance $s \in \RR$ from the origin $\bo$. Observe that this notation is justified by rotational symmetry around $\bo$. In \cite[Proposition 4.1]{KRT} it is proven that this intersection is empty for $|s| \geq R$, and otherwise satisfies
	\begin{equation}\label{eq:vol_intersection_horo}
		\cH^{d-1}(H(s) \cap B_R^d) = \kappa_{d-1} \bigl[ 2 e^s (\cosh R - \cosh s)\bigr]^{\frac{d-1}{2}},
	\end{equation}
	where for an integer $\ell\geq 1$ we write $\kappa_\ell$ for the volume of the $\ell$-dimensional Euclidean unit ball.
	
	{We also mention some elementary properties of the Kolmogorov and Wasserstein metric that will be useful for us. First, for any random variables $X$ and $Y$, and scalars $\alpha>0$ one has
	\begin{equation}\label{eq:Wass_homogeneity}
	\begin{aligned}
	d_{\Wass}(\alpha X, \alpha Y) &= \alpha \, d_{\Wass}(X,Y),\\
	d_{\Kol}(\alpha X, \alpha Y) &= d_{\Kol}(X,Y).
	\end{aligned}
	\end{equation}
	Moreover, given independent random variables $X$ and $Y$ as well as another random variable $Z$ with density bounded by $1$, one has
	\begin{equation}\label{eq:Wass_sums}
	d_{\bullet}(X+Y,Z) \leq d_{\bullet}(X,Z) + \EE |Y|.
	\end{equation}	
	For the Wasserstein metric this follows at once from the triangle inequality and the fact that the Wasserstein metric is bounded by the $L^1$-metric (and, indeed, holds without the independence and bounded density assumption). For the Kolmogorov metric this is a consequence of the following result.
	
	\begin{lemma}\label{lem:Kol_sum}
	Let $X$ and $Y$ be two independent random variables, and let $Z$ be another random variable admitting a bounded density $f_Z$. Then
	\[
	d_{\Kol}(X+Y,Z) \leq d_{\Kol}(X,Z) + \|f_Z\|_\infty \cdot \EE |Y|.
	\]
	\end{lemma}}
	
	{\begin{proof}
	We write $F_W$ for the cumulative distribution function of the random variable $W$. For any $t\in\RR$ one has, using independence in the second step,
	\begin{align*}
	|\PP(X+Y\leq t)-\PP(Z\leq t)| &= |\EE\left[\PP(X\leq t-Y|Y)\right] - \PP(Z\leq t)| \\
	&\leq |\EE\left[F_X(t-Y)-F_Z(t-Y)\right] | + |\EE[F_Z(t)-F_Z(t-Y)]| \\ & \leq d_{\Kol}(X,Z) + \|f_Z\|_\infty\cdot \EE |Y|.
	\end{align*}
	Taking the supremum over $t\in\RR$ yields the  result.
	\end{proof}
	\begin{remark}
	As is evident from the proof, the assumption that $Z$ has density may be removed, in which case the  result takes the  form
		\[
	d_{\Kol}(X+Y,Z) \leq d_{\Kol}(X,Z) +  \EE [\omega_Z(|Y|)],
	\]
	where $\omega_Z$ denotes the modulus of continuity of $F_Z$, that is, $\omega_Z(\eps) := \sup_{t\in\RR} \PP(t< Z\leq t+\eps)$.
	\end{remark}
}
	

	{A first step in the proof of our main result} is to reduce the normal approximation bound to the following integral estimate. Define
	\[
	J_{R,d} := \int_0^R \left(1 - {\cosh s - 1 \over \cosh R - 1}\right)^{d-1}\,\dint s.
	\]
	\begin{proposition}\label{prop:Wass_bd_J}
	The following bound holds true for all $d\geq 2$ and $R>0$ {and any $\bullet \in \{\Kol, \Wass\}$:
\begin{equation}\label{eq:bd_J_simple}
d_{\bullet}\left( \frac{S_{R,d} - \EE S_{R,d}  } {\sqrt{\Var S_{R,d} } }, N_{\tinyonehalf} \right)   \leq {c} \cdot \left( {1 \over \sqrt{d}\, J_{R,d}}  + {1 \over d\,\sqrt{J_{R,d}} }\right).
\end{equation}

	for some universal constant $c\in(0,\infty)$.}

	\end{proposition}

	\begin{proof}
				In view of the representation \eqref{eq:dHlambda} of the invariant measure $\Lambda$ and the expression \eqref{eq:vol_intersection_horo} for the intersection volume, we have that
		\begin{equation}\label{eq:S-sum-dist-process}
			S_{R,d} = \sum_{s \in \xi} f_R(s),
		\end{equation}
		where $\xi$ is an inhomogeneous Poisson process on $\RR$ with density $s \mapsto e^{-(d-1)s}$, and the function $f_R$ is defined by 
		\begin{equation}\label{eq:fR}
		f_R(s) = \begin{cases}
		\kappa_{d-1} \bigl[ 2 e^s (\cosh R - \cosh s)\bigr]^{\frac{d-1}{2}}  &: |s| \leq R, \\ 0 &: \text{else.}
		\end{cases}
		\end{equation}
		We decompose the random variable $S_{R,d}$ into a `positive' and `negative' part as follows:
		\begin{equation*}
			S_{R,d} = S_{R,d}^{+} + S_{R,d}^{-},
		\end{equation*}
		where
		\begin{equation*}
			S_{R,d}^{+} := \sum_{\substack{s \in \xi \\  s > 0}} f_R(s) \qquad \text{and} \qquad S_{R,d}^{-} := \sum_{\substack{s \in \xi \\  s < 0}} f_R(s).
		\end{equation*}
		{Note that from the independence property of Poisson processes it follows that the random variables $S_{R,d}^{+}$ and $	S_{R,d}^{-}$ are independent.} We then have
		\begin{equation*}
			\frac{S_{R,d} - \EE S_{R,d} }{\sqrt{\Var S_{R,d} }} = \frac{S_{R,d}^{+} - \EE S_{R,d}^{+}}{\sqrt{\Var S_{R,d}}} + \frac{S_{R,d}^{-} - \EE S_{R,d}^{-}}{\sqrt{\Var S_{R,d}}}.
		\end{equation*}
		Observe that 
		\begin{equation}\label{eq:vardecomp}
			\Var S_{R,d}^{+} = \Var S_{R,d}^{-} = {1 \over 2} \, \Var S_{R,d},
		\end{equation}
		which follows from the evenness of the integrand in the variance representation
		\begin{equation*}
			\Var S_{R,d} = \int_\RR f_R^2(s) e^{-(d-1)s}\,\dint s  = 2^{d-1} \kappa_{d-1}^2 \int_{-R}^R (\cosh R - \cosh s)^{d-1} \,\dint s,
		\end{equation*}
		which in turn follows from the multivariate Mecke formula for Poisson processes \cite[Theorem 4.1]{LP}, where we used \eqref{eq:dHlambda} and \eqref{eq:vol_intersection_horo}.	 We deduce that
		\begin{equation}\label{eq:bound_sum}
		{d_{\bullet}\left( \frac{S_{R,d} - \EE S_{R,d}  } {\sqrt{\Var S_{R,d} } }, N_{\tinyonehalf} \right)  \leq 
		 d_{\bullet}\left( \frac{S_{R,d}^{-} - \EE S _{R,d}^{-}  } {\sqrt{\Var S_{R,d}^-  } }, N \right) + 2^{-{1\over 2}}\, \EE \left|\frac{S_{R,d}^{+} - \EE S _{R,d}^{+}  } {\sqrt{\Var S_{R,d}^+  } } \right|,}
		\end{equation}
		where $N$ is a standard Gaussian random variable, and where we have used \eqref{eq:Wass_homogeneity} and \eqref{eq:Wass_sums} together with the fact that $2^{-{1\over 2}} N $ has the same distribution as our target random variable $N_{\tinyonehalf}$.
		
		To control the first summand in \eqref{eq:bound_sum}, we apply {the following normal approximation bound, which is a special case of general bounds for so-called Poisson $U$-statistics, see {\cite[Theorem 3.12 and Equation (3.9)]{LrP} and \cite[page 112]{Schulte}}}. Applied to $S _{R,d}^{-}$, it states that
		{\begin{equation*}
		d_{\bullet }\left(\frac{S _{R,d}^{-}- \EE S _{R,d}^{-} }{\sqrt{\Var S _{R,d}^{-} } },N\right) \leq c_\bullet {{ \sqrt{\cum_4(S_{R,d}^-)}} \over \Var S_{R,d}^- },
		\end{equation*}
		for some constants $c_\bullet \in(0,\infty)$, $\bullet \in \{\Kol,\Wass \}$ (explicitly, one can take $c_{\Kol}=19$ and $c_{\Wass}=2$)}, where {$\cum_4(W) = \EE\left[(W-\EE W)^4\right]-3\Var(W)^2$ denotes the fourth cumulant of a random variable $W$}. 
		{Noting that $S_{R,d}^+ \geq 0$, the} second summand in \eqref{eq:bound_sum} is easily bounded by 
			\begin{align*}
					\frac{  \EE| S_{R,d}^{+} - \EE S _{R,d}^{+}  | }{ \sqrt{\Var S_{R,d}^+  } } &
					\leq {2\,\EE S_{R,d}^+\over  \sqrt{\Var S_{R,d}^+  } }.
				\end{align*}
		This gives
		\begin{align}
		d_{\bullet}\left( \frac{S_{R,d} - \EE S_{R,d}  } {\sqrt{\Var S_{R,d} } }, N_{\tinyonehalf} \right) 
		& \leq {c_{\bullet}} {\sqrt{\cum_4(S_{{R,d}}^-)} \over \Var S_{{R,d}}^-} + 2^{1/2} {\EE S_{{R,d}}^+ \over \sqrt{\Var S_{{R,d}}^+}}.\label{eq:Wass_bd_cum}
		\end{align}
		If we denote further $C_d: = 2^{(d-1)/2} \kappa_{d-1}$ and define
		\begin{align*}
		I_1(R) &:= \int_0^R (\cosh R - \cosh s)^{d-1 \over 2}e^{-{d-1\over 2}s}\,\dint s,\\
		I_2(R) &:= \int_0^R (\cosh R - \cosh s)^{d-1}\,\dint s,\\
		I_3(R) &:= \int_0^R (\cosh R - \cosh s)^{2(d-1) }e^{-{(d-1)}s}\,\dint s,
		\end{align*}
		then we compute, using \eqref{eq:S-sum-dist-process} and \eqref{eq:fR}, that
		\begin{align*}
		\EE S_{R,d}^+ = C_d I_1(R), \qquad \Var(S_{R,d}^\pm) = {C_d^2} I_2(R), \qquad \cum_4 (S_{R,d}^-) = {C_d^4} I_3(R).
		\end{align*}
		Here the expectation and variance are computed with the help of the multivariate Mecke equation, and the fourth cumulant using \cite[Corollary 1]{LPST}. Plugging this into \eqref{eq:Wass_bd_cum} finally gives
		\begin{equation}\label{eq:Wass_bd_Is}
	{d_{\bullet}}\left( \frac{S_{R,d} - \EE S_{R,d}  } {\sqrt{\Var S_{R,d} } }, N_{\tinyonehalf} \right)  
		\leq {c_\bullet} \left[{\sqrt{I_3(R)}\over I_2(R)} + {I_1(R) \over \sqrt{I_2(R)}}\right]
		\end{equation}
		(since in both cases, $c_\bullet > 2^{1/2}$.) Now we use the following trivial estimates for $I_1$ and $I_3$:
		\begin{align*}
		I_1(R) & \leq (\cosh R - 1)^{d-1 \over 2} \cdot {2 \over d-1}, \\
		I_3(R) & \leq (\cosh R-1)^{2(d-1)} \cdot {1 \over d-1}.
		\end{align*}
		Moreover, for $I_2$ we write
		\begin{align*}
		I_{2}(R) &= \int_0^R (\cosh R - \cosh s)^{d-1} \,\dint s \\
		&= (\cosh R - 1)^{d-1} \int_0^R \left(1 - {\cosh s - 1 \over \cosh R - 1}\right)^{d-1}\,\dint s \\ 
		& = (\cosh R - 1)^{d-1} J_{R,d}.
		\end{align*}
		Plugging all this back into \eqref{eq:Wass_bd_Is} leads to the estimate
		\begin{equation*}
		{d_{\bullet}}\left( \frac{S_{R,d} - \EE S_{R,d}  } {\sqrt{\Var S_{R,d} } }, N_{\tinyonehalf} \right)   \leq {c}_\bullet \left( {1 \over \sqrt{d-1}\, J_{R,d}}  + {2 \over (d-1)\,\sqrt{J_{R,d}} }\right),
		\end{equation*}
		which clearly implies the bound \eqref{eq:bd_J_simple} and completes the proof.
	\end{proof}
		
		Our next task therefore is to estimate $J_{R,d}$. This is achieved by the following result.
		
		{\begin{lemma}\label{lem:bound-J-general}
		There exists a constant $C>0$ such that the following holds for any $d \geq 2$ and $R\geq 1$.
		\begin{equation*}
		J_{R,d} \geq C\cdot \begin{cases}
		{e^{R/2} \over \sqrt{d}}  &: R-\log d \leq 1,\\
		R-\log d &: R-\log d >1.
		\end{cases}
		\end{equation*}	
		\end{lemma}}
	
		{\begin{remark}
		As in Theorem \ref{thm:CLT-master} (see Remark \ref{rem:main_thm}), the threshold $1$ between the two cases of Lemma \ref{lem:bound-J-general} can be replaced by any positive number, in which case the constant $C$ may change as well.
		\end{remark}}
	
		We postpone the proof of Lemma \ref{lem:bound-J-general} until Section \ref{sec:J}, and first use it to deduce our main result.
	
	\begin{proof}[Proof of Theorem \ref{thm:CLT-master}]
		The theorem follows upon combining Proposition \ref{prop:Wass_bd_J} with the integral estimate provided by Lemma \ref{lem:bound-J-general}. {Indeed, in the case $R-\log d >1$, we obtain immediately the bound asserted by the theorem. In the case $R-\log d \leq 1$, combining Proposition \ref{prop:Wass_bd_J} with the first case of Lemma \ref{lem:bound-J-general} gives
		\[
		{d_{\bullet}}\left( \frac{S_{R} - \EE S_{R}  } {\sqrt{\Var S_{R} } }, N_{\tinyonehalf} \right)   \leq C \left( e^{-R/2} + d^{-3/4} e^{-R/4}\right).
		\]
		Since by assumption, $e^{R} \leq e \cdot d$, the first term {on} the right-hand-side is dominant, leading to the bound appearing in the corresponding case of the theorem. }
	\end{proof}

\section{Bounding $J_{{R,d}}$}\label{sec:J}

Here we bound the integral $ J_{R,d}$, which we recall is given by
\[
J_{R,d} = \int_0^R \left(1 - {\cosh s - 1 \over \cosh R - 1}\right)^{d-1}\,\dint s.
\]

\begin{proof}[Proof of Lemma \ref{lem:bound-J-general}]
First we write, with the help of the hyperbolic identity $\cosh x -1 = 2 \sinh^2 {x \over 2}$,
\begin{align}
	J_{R,d} & = \int_0^{R} \left(1 - {\sinh^2(s/2) \over \sinh^2(R/2)}\right)^{d-1}\,\dint s. \label{eq:Jd_sinh}
\end{align}
	Next we define  $$ \rho  = \rho(R,d):= {\sinh(R/2) \over \sqrt{d}}.$$
	We now proceed with the two cases of the lemma separately.
	
	{\begin{itemize}
	\item[1.] Suppose first that $R-\log d \leq 1$. In the integral \eqref{eq:Jd_sinh} we make the substitution 
	\[
	s = 2 \asinh (\rho x); \quad \dint s = {2\rho \over \sqrt{1+\rho^2x^2}} \,\dint x,
	\]
	so that 
	\[
	{\sinh (s/2) \over \sinh(R/2)} = {x \over \sqrt d}.
	\]	
	 Note that the assumption $R-\log d \leq 1$ implies 
	\[
	\rho \leq {e^{R/2} \over 2 \sqrt{d}} = {1 \over 2} \exp\left({R-\log d \over 2}\right) \leq 1,
	\]
	which altogether yields
	\begin{align*}
	J_{R,d} &= 2\rho \int_0^{\sqrt{d}} \left(1-{x^2\over d}\right)^{d-1} \,{\dint x \over {\sqrt{1+\rho^2 x^2}}} \\ &\geq 2\rho \int_0^{\sqrt{d}} \left(1-{x^2\over d}\right)^{d-1} \,{\dint x \over {\sqrt{1+ x^2}}}
	\end{align*}
	Finally, noting that the latter sequence of integrals converges as $d\to\infty$ to 
	\[
	\int_0^\infty e^{-x^2} \,{\dint x\over \sqrt{1+x^2}} >0,
	\]
	we deduce that $J_{R,d} \geq \hat{C} \cdot \rho = \hat{C} \cdot {\sinh(R/2)\over \sqrt{d}}$ for some constant $\hat{C}>0$. As $R \geq 1$, this implies that $J_{R,d} \geq C \cdot {e^{R/2} \over \sqrt{d}}$ for some $C > 0$, as asserted.
	
	\item[2.] Consider now the case $R-\log d >1$. Noting that, by definition of $\rho$, ${s \leq 2 \asinh (\rho)}$ implies ${\sinh(s/2)\over \sinh(R/2)} \leq {1\over \sqrt d}$, we estimates \eqref{eq:Jd_sinh} by
	\begin{align*}
		J_{R,d} & \geq \int_0^{2\,\asinh (\rho)} \left(1 - {1 \over d}\right)^{d-1}\,\dint s\geq {2 \over e} \asinh (\rho).
	\end{align*}

Now, using the logarithmic representation of the inverse hyperbolic sine
	\[
	\asinh (x) = \log\left(x+\sqrt{x^2+1}\right) \geq \log(2x),
	\]
	we deduce that
	\begin{align*}
		\asinh(\rho) \geq \log \left({e^{R/2}-e^{-R/2} \over \sqrt{d}}\right) = {R-\log d \over 2}+\log(1-e^{-R}).
	\end{align*}
	Finally, the assumption $R-\log d \geq 1$ gives 
	\[
	\log(1-e^{-R}) \geq \log\left(1-{1\over 2e}\right) \geq -{R-\log d \over 4},
	\]
	which implies that
	\[
	\asinh (\rho) \geq {R-\log d \over 4}.
	\]
	This shows that $J_{R,d} \geq C (R-\log d)$, and completes the proof.\hfill \qedhere
	\end{itemize}}
\end{proof}

\section{The Euclidean case}\label{sec:Euclidean}

Let $\eta_d$ be a stationary and isotropic Poisson process on the space $\AA(d,d-1)$ of affine hyperplanes in $\RR^d$ with intensity $1$. Its intensity measure $\Lambda_e$ is then given by
$$
\int_{\AA(d,d-1)}f(H)\,\Lambda_e(\dint H) = \int_{\RR}\int_{\SS^{d-1}}f(H_e(s,u))\,\dint u\,\dint s,
$$
for a non-negative measurable function $f:\AA(d,d-1)\to\RR$, where $H_e(s,u)$ stands for the unique hyperplane in $\RR^d$ with signed distance $s$ from $\bf o$ and unit normal vector $u$. {{As above, $\dint s$ and $\dint u$ stand for the Lebesgue measure on $\RR$ and the normalized spherical Lebesgue measure on $\SS^{d-1}$, respectively}}. By 
$$
S_{R,d,e} := \sum_{H\in\eta_d}{\cH_e^{d-1}}(H\cap B_{R,e}^d)
$$
we denote the total surface area induced by the hyperplanes of $\eta_d$ within a centred Euclidean ball $B_{R,e}^d$ of radius $R>0$, where {the Hausdorff measure $\cH_e^{d-1}$ is understood with respect to the Euclidean metric}. Using, as above \cite[Theorem 3.12 and Equation (3.9)]{LrP} and \cite[page 112]{Schulte}
we find that
\begin{align}\label{eq:EucliBound}
{d_{\bullet}}\Big({S_{R,d,e}-\EE S_{R,d,e}\over\sqrt{\Var S_{R,d,e}}},N\Big) \leq {c_{\bullet}}{\sqrt{\cum_4(S_{R,d,e})}\over\Var S_{R,d,e}}
\end{align}
with a standard Gaussian random variable $N$. The variance and the fourth cumulant of $S_{R,d,e}$ are given explicitly by
\begin{align*}
\Var S_{R,d,e} &= \int_{\cH_e}{\cH_e^{d-1}}(H\cap B_{R,e}^d)^2\,\Lambda_e(\dint H),\\
\cum_4(S_{R,d,e}) &= \int_{\cH_e}{\cH_e^{d-1}}(H\cap B_{R,e}^d)^4\,\Lambda_e(\dint H).
\end{align*}
Denoting by $\kappa_{d-1}$ the $(d-1)$-volume of the $(d-1)$-dimensional Euclidean unit ball, we have that
\begin{align*}
\Var S_{R,d,e} &= 2\kappa_{d-1}^2\int_0^R(R^2-s^2)^{d-1}\,\dint s\\
&=2\kappa_{d-1}^2R^{2d-1}\int_0^1(1-t^2)^{d-1}\,\dint t\\
&={\pi^{d-{1\over 2}}\Gamma(d)R^{2d-1}\over\Gamma({d\over 2}+{1\over 2})^2\Gamma(d+{1\over 2})},
\end{align*}
where we applied the substitution {$s\mapsto Rt$}. The same computation also leads to an explicit expression for $\cum_4(S_{R,d,e})$:
\begin{align*}
\cum_4(S_{R,d,e}) &= 2\kappa_{d-1}^4\int_0^R(R^2-s^2)^{2(d-1)}\,\dint s\\
&=2\kappa_{d-1}^4R^{4d-3}\int_0^1(1-t^2)^{2(d-1)}\,\dint t\\
&={\pi^{2d-{3\over 2}}\Gamma(2d-1)R^{4d-3}\over\Gamma({d\over 2}+{1\over 2})^4\Gamma(2d-{1\over 2})}{.}
\end{align*}
In conjunction with \eqref{eq:EucliBound} this gives 
\begin{align*}
d_{\rm \bullet}\Big({S_{R,d,e}-\EE S_{R,d,e}\over\sqrt{\Var S_{R,d,e}}},N\Big) \leq 
{2^{-1/4} c_{\bullet}}{\Gamma(d+{1\over 2})\over\Gamma(d)}\sqrt{{\Gamma(2d-1)\over\Gamma(2d-{1\over 2})}}\,R^{-1/2}.
\end{align*}
Using the well-known asymptotics for quotients of gamma functions, as $d\to\infty$ we arrive at
$$
d_{\rm \bullet}\Big({S_{R,d,e}-\EE S_{R,d,e}\over\sqrt{\Var S_{R,d,e}}},N\Big) \leq C\,d^{1/4}R^{-1/2}
$$
for some absolute constant $C>0$.

\subsection*{Acknowledgement}

We wish to thank Matthias Schulte (Hamburg) for motivating us to study the problem addressed in this paper. {We also thank two anonymous referees for their valuable comments that helped us to improve our text.}\\ DR and CT were supported by the German Research Foundation (DFG) via CRC/TRR 191 \textit{Symplectic Structures in Geometry, Algebra and Dynamics}. ZK and CT were supported by the German Research Foundation (DFG) via the Priority Program SPP 2265 \textit{Random Geometric Systems}. ZK was also supported by the German Research Foundation (DFG) under Germany’s Excellence Strategy EXC 2044 – 390685587 \textit{Mathematics M\"unster: Dynamics -- Geometry -- Structure}.

\end{document}